\documentstyle[12pt,amssymb,amstex]{amsart}

\numberwithin{equation}{section}

\def\ce{{\cal E}}

\def\cs{{\cal S}}

\def\ga{{\frak A}}
\def\gb{{\frak B}}

\def\bc{{\mathbb C}}
\def\bbf{{\mathbb F}}

\def\bd{{\mathbb D}}

\def\bn{{\mathbb N}}

\def\br{{\mathbb R}}
\def\bt{{\mathbb T}}

\def\bz{{\mathbb Z}}

\def\a{\alpha}
\def\b{\beta}
  
\def\d{\delta}  

\def\e{\epsilon}
\def\l{\lambda} 

\def\m{\mu}
\def\p{\psi}
\def\n{\nu}

\def\s{\sigma} 
\def\t{\tau}
\def\f{\varphi} 
\def\v{\phi}
  
\def\w{\omega} \def\Om{\Omega}

\newtheorem{thm}{Theorem}[section]

\newtheorem{cor}[thm]{Corollary}
\newtheorem{prop}[thm]{Proposition}

\theoremstyle{remark}
\newtheorem{rem}{Remark}[section]

\def\di{\mathop{\rm d}}

\def\id{\mathop{\rm id}}

\def\id{{\bf 1}\!\!{\rm I}}

\begin{document}

\small

\title[On strictly weak mixing]
{On strictly weak mixing $C^*$-dynamical systems and\\ a weighted
ergodic theorem}

\thanks{One leave from: National University of Uzbekistan, Dep. Mechanics \& Mathematics, Tashkent, Uzbekistan}

\author{Farrukh Mukhamedov}
\address{Farrukh Mukhamedov\\
 Department of Computational \& Theoretical Sciences\\
Faculty of Sciences, International Islamic University Malaysia\\
P.O. Box, 141, 25710, Kuantan\\
Pahang, Malaysia} \email{{\tt far75m@@yandex.ru}}

\begin{abstract}
We prove that unique ergodicity of tensor product of
$C^*$-dynamical system implies its strictly weak mixing. By means
of this result a uniform weighted ergodic theorem with respect to
$S$-Besicovitch sequences for strictly weak mixing dynamical
systems is proved. Moreover, we provide certain examples of
strictly weak mixing dynamical systems.

\vskip 0.3cm \noindent {\it Mathematics Subject Classification}:
46L35, 46L55, 46L51, 28D05
60J99.\\
{\it Key words}: unique ergodicity, strict weak mixing,
$C^*$-dynamical system, $S$-Besicovitch sequence.
\end{abstract}

\maketitle
\section{Introduction}

Recently, the investigation of the ergodic properties of quantum
dynamical systems had a considerable growth. Since the theory of
quantum dynamical systems provides convenient mathematical
description of the irreversible dynamics of an open quantum system
(see \cite{BR}, sec.4.3, \cite{Wa},\cite{Ru}). In this setting,
the matter is more complicated than in the classical case. Some
differences between classical and quantum situations are pointed
out in \cite{AH},\cite{NSZ}. This motivates an interest to study
of dynamics of quantum systems (see
\cite{AH},\cite{FR1},\cite{J}). Therefore, it is then natural to
address the study of the possible generalizations to quantum case
of the various ergodic properties known for classical dynamical
systems.  A lot of papers (see, \cite{FV},
\cite{FR2},\cite{L1},\cite{L2},\cite{W}) were devoted to the
investigations of mixing properties of dynamical systems.

It is known \cite{KSF} that a strong ergodic property for a
classical system is the unique ergodicity. Namely, a classical
dynamical system $(\Om,T)$ consisting of  a compact Hausdorff
space $\Om$ and a homeomorphism $T$ is said to be {\it uniquely
ergodic} if there exists a unique invariant Borel measure $\m$ for
$T$. It is seen that the ergodic average
${\displaystyle\frac1n\sum_{k=0}^{n-1}f\circ T^{k}}$ converges
uniformly to the constant function $\int f\di\m$ in this case. A
pivotal example of classical uniquely ergodic dynamical system is
given by an irrational rotation on the unit circle, see e.g.
\cite{KSF}. In quantum setting, the last property is formulated as
follows (see also Sec. 2). Let $(\ga,T)$ be a $C^{*}$--dynamical
system based on the $C^{*}$--algebra $\ga$ and a unital completely
positive (ucp) map $T$ on $\ga$. The unique ergodicity or
equivalently strict ergodicity for $(\ga,\a)$ is equivalent (cf.
\cite{AD, MT}) to the norm convergence
\begin{equation}
\label{er}
\lim_{n\to+\infty}\frac1n\sum_{k=0}^{n-1}T^{n}(a)=E(a)\,,\quad
a\in\ga\,,
\end{equation}
where $E$ is a  conditional expectation, given by
$E=\f(\,\cdot\,)\id$, onto the fixed--point subspace of $T$,
consisting of the constant multiples of the identity. Here,
$\f\in\cs(\ga)$ is the unique invariant state for $T$. Some
generalizations of unique ergodicity have been investigated in
\cite{AD,AM}, where the conditional expectation $E$ in \eqref{er}
(necessarily unique) is taken as a projection onto the
fixed--point subspace of $T$, which, in general, is supposed to be
nontrivial.

In \cite{MT} we have introduced a property stronger than the
unique ergodicity, called {\it strict weak mixing}. This property
for $(\ga,T)$ -$C^{*}$--dynamical system requires the existence of
a state $\f\in\cs(\ga)$ such that
\begin{equation}
\label{mta}
\lim_{n}\frac1n\sum_{k=0}^{n-1}\big|\psi(T^{k}(a))-\f(a)\big|=0\,,
\quad a\in\ga\,,
\end{equation}
for each $\psi\in\cs(\ga)$. It can be shown (see below) that $\f$
is the unique invariant state for $T$. If $(\ga,T)$ is strictly
weak mixing, then it is uniquely ergodic (see Proposition
\ref{uesm}) . Conversely, the irrational rotations on the unit
circle provide examples of uniquely ergodic dynamical systems
which are not strictly weak mixing, see \cite{MT}, Example 2.

In this paper we are going to prove for the strict mixing an
analog of the well-known classical result stating that a
transformation is weakly mixing if and only if its Cartesian
square is ergodic \cite{ALW}.

The paper is organized as follows. In section 2 we give a needed
preliminary definitions and results. In section 3 we will prove
that if tensor product of dynamical system is uniquely ergodic,
then it is strictly weakly mixing. In section 4, we provide
certain examples of strictly weak mixing dynamical systems. In the
final Section 5, by means of the main result we prove one uniform
weighted ergodic theorem for strictly weak mixing systems.

\section{Preliminaries}

In this section we recall some preliminaries concerning
$C^*$-dynamical systems.

Let $\ga$ be a $C^*$-algebra with unit $\id$. An element $x\in\ga$
is called {\it self-adjoint} (resp. {\it positive}) if $x=x^*$
(resp. there is an element $y\in\ga$ such that $x=y^*y$). The set
of all self-adjoint (resp. positive) elements will be denoted by
$\ga_{sa}$ (resp. $\ga_+$). By $\ga^*$ we denote the conjugate
space to $\ga$. A linear functional $\f\in\ga^*$ is called {\it
Hermitian} if $\f(x^*)=\overline{\f(x)}$ for every $x\in\ga$. A
Hermitian functional $\f$ is called {\it positive} if
$\f(x^*x)\geq 0$ for every $x\in\ga$. A positive functional $\f$
is said to be a {\it state} if $\f(\id)=1$. By $\cs(\ga)$ (resp.
$\ga^*_h$) we denote the set of all states (resp. Hermitian
functionals) on $\ga$. Let $\gb$ be another $C^*$-algebra with
unit. By $\ga\odot\gb$ we denote the algebraic tensor product of
$\ga$ and $\gb$. A completion of $\ga\odot\gb$ with respect to the
minimal $C^*$-tensor norm on $\ga\odot\gb$ is denoted by
$\ga\otimes\gb$, and it would be also a $C^*$-algebra with a unit
(see, \cite{T}). A linear operator $T:\ga\mapsto\ga$ is called
{\it positive} if $Tx\geq 0$ whenever $x\geq 0$. By $M_n(\ga)$ we
denote the set of all $n\times n$-matrices $a=(a_{ij})$ with
entries $a_{ij}$ in $\ga$. A linear mapping $T:\ga\mapsto\ga$ is
called {\it  completely positive} if the linear operator
$T_n:M_n(\ga)\mapsto M_n(\ga)$ given by $T_n(a_{ij})=(T(a_{ij}))$
is positive for all $n\in\bn$. A completely positive map
$T:\ga\mapsto\ga$ with $T\id=\id$ is called a {\it unital
completely positive (ucp)} map.  A pair $(\ga,T)$ consisting of a
$C^*$-algebra $\ga$ and a ucp map $T:\ga\mapsto\ga$ is called {\it
a $C^*$-dynamical system}. In the sequel, we will call any triplet
$(\ga,\f,T)$ consisting of a $C^*$-algebra $\ga$, a state $\f$ on
$\ga$ and a ucp map $T:\ga\mapsto\ga$ with $\f\circ T=\f$, that is
a dynamical system with an invariant state, {\it a state
preserving $C^*$-dynamical system}. A state preserving
$C^*$-dynamical system is a non-commutative $C^*$-probability
space $(\ga,\f)$ (see \cite{CO}) together with a ucp map $T$ on
$\ga$ preserving the non-commutative probability $\f$. It is known
\cite{T} that if $(\ga,T)$ and $(\gb,H)$ are two  $C^*$-dynamical
systems, then $(\ga\otimes\gb,T\otimes
H)$ is also $C^*$-dynamical system. Since a mapping $T\otimes
H:\ga\otimes\gb\mapsto\ga\otimes\gb$ given by $(T\otimes
H)(x\otimes y)=Tx\otimes Hy$ is a ucp map.

We say that the state preserving $C^*$-dynamical system
$(\ga,\f,T)$ is {\it ergodic} (respectively, {\it weakly mixing,
strictly weak mixing}) with respect to $\f$ if
\begin{equation}\label{erg}
\lim_{n\to\infty}\frac{1}{n}\sum_{k=0}^{n-1}(\f(yT^k(x))-\f(y)\f(x))=0,
\ \ \textrm{for all} \ \ x,y\in\ga. \end{equation}

(respectively, \begin{equation}\label{wmix}
\lim_{n\to\infty}\frac{1}{n}\sum_{k=0}^{n-1}|\f(yT^k(x))-\f(y)\f(x)|=0,\
\ \textrm{for all} \ \ x,y\in\ga, \end{equation}

\begin{equation}\label{stmix1}
\lim_{n\to\infty}\frac{1}{n}\sum_{k=0}^{n-1}|\p(T^k(x))-\p(\id)\f(x)|=0,
\ \ \textrm{for all} \ \ x\in\ga, \p\in\ga^*.) \end{equation}

The  state preserving $C^*$-dynamical system $(\ga,\f,T)$ is
called {\it uniquely ergodic} with respect to $\f$ if $\f$ is the
unique invariant state under $T$.

\begin{rem} If we take a functional $\f(xy)$ instead of
$\p(x)$ in \eqref{stmix1}, then one can see that strict weak
mixing implies weak mixing. Converse, is not true. A related
example was provided in \cite{MT}, Example 3.
\end{rem}

In \cite{MT} (see also \cite{AD}) we have proved the following
characterization of unique ergodicity of dynamical systems.

\begin{thm}\label{ue}  Let $(\ga,\f,T)$ be a state preserving
$C^*$-dynamical system. The following conditions are equivalent
\begin{itemize}
\item[(i)] $(\ga,\f,T)$ is uniquely ergodic ; \item[(ii)] For
every $x\in\ga$ the following equality holds
$$
\lim_{n\to\infty}\frac{1}{n}\sum_{k=0}^{n-1}T^k(x)=\f(x)\id,
$$
where convergence in norm of $\ga$; \item[(iii)] For every
$x\in\ga$ and $\p\in \ga^*$ the following equality holds
$$
\lim_{n\to\infty}\frac{1}{n}\sum_{k=0}^{n-1}\p(T^k(x))=\p(\id)\f(x).
$$
\end{itemize}
\end{thm}

\begin{rem} From this Theorem we immediately infer that
unique ergodicity implies ergodicity of $C^*$-dynamical system.
\end{rem}

\begin{prop}\label{uesm} If the  $C^*$-dynamical system $(\ga,\f,T)$ is
strictly weak mixing, then it is uniquely ergodic.
\end{prop}
\begin{pf}
Let $\p\in \ga^*$, then one gets
$$
\bigg|\frac{1}{n}\sum_{k=0}^{n-1}\big(\p(T^{k}(x))-\p(\id)\f(x)\big)\bigg|
\leq\frac{1}{n}\sum_{k=0}^{n-1}\big|\p(T^{k}(x))-\p(\id)\f(x)\big|\to0
$$
whenever $n\to\infty$, as $(\ga,T)$ is strictly weak mixing. By
using the Jordan decomposition of bounded linear functionals (cf.
\cite{T}), we conclude that (iii) of Theorem \ref{ue} is
satisfied.
\end{pf}

In many interesting situations, the ergodic behavior of dynamical
systems is connected with some spectral properties, see e.g.
\cite{C, L, NSZ,W}. It is not possible to extend such results to
our situation. However, a strictly weak mixing map $T$ cannot have
eigenvalues on the unit circle $\bt$ except $z=1$.

Let $\bd:=\{z\in\bc\,:\,|z|\leq 1\}$ be the unit disk in the
complex plane, and $\stackrel{\circ}{\bd}=\{z\in\bc\,:\,|z|<1\}$
its interior. If $T$ has norm one, we have $\s(T)\subset\bd$,
$\s(T)$ being the spectrum of $T$.

Let $T:\ga\mapsto\ga$ be a linear map. Denote
\begin{eqnarray*}
&&\ga_z=\{x\in\ga: T(x)=zx\},\\
&& \ga_z^*=\{f\in\ga^*: f\circ T=z f \},
\end{eqnarray*}
where $z\in\bc$. Furthermore,
\begin{prop}
\label{mix2} Let $\big(\ga,T\big)$ be a strictly weak mixing
$C^*$-dynamical system. Then $z\in\bt\backslash\{1\}$ implies
$\ga_z=\{0\}$ and $\ga_z^*=\{0\}$
\end{prop}
\begin{pf}
Assume that $T(x_0)=zx_0$ for some $z\neq 1$. Then
$\f(x_0)=\f(T(x_0))=z\f(x_0)$ which means $\f(x_0)=0$. In
addition, the strict weak mixing implies
\begin{align*}
0=&\lim_{n}\frac{1}{n}\sum_{k=0}^{n-1}\big|\p(T^k(x_0))-\p(\id)\f(x_0)\big|
=\lim_{n}\frac{1}{n}\sum_{k=0}^{n-1}\big|z^{k}\p(x_0)\big|\\
=&\lim_{n}\frac{1}{n}\sum_{k=0}^{n-1}\big|\p(x_0)\big| =
|\p(x_0)|\,.
\end{align*}
Namely, $\p(x_0)=0$ for every $\p\in \ga^*$, hence $x_0=0$. The
second part can be proceeded similarly.
\end{pf}

\begin{rem} For any linear map $T$ of $\ga$, it is obvious that
if $z\in \stackrel{\circ}{\bd}$ and $x\in\ga_z$, then
${\displaystyle\lim_{k}T^{k}(x)=0}$.
\end{rem}

\section{Tensor product of strictly weak mixing dynamical systems}

This section is devoted to tensor product of uniquely ergodic and
strictly weak mixing dynamical systems. Here we prove the main
result of the paper.

Set
$$
\ga^*_1=\{g\in\ga^*: \ \|g\|_1\leq 1\}, \ \ \
\ga^*_{1,h}=\ga^*_1\cap\ga^*_h.
$$

Now we are going to prove an analogous result of \cite{ALW,Wa} for
the strictly weak mixing dynamical systems.

\begin{thm}
\label{mix-a} Let $(\ga,\f,T)$, $(\gb,\f_1,H)$ be two state
preserving $C^*$-dynamical systems.   For the following assertions
\begin{itemize}
\item[(i)] The state preserving  $C^*$-dynamical system
$(\ga\otimes\gb,\f\otimes\f_1,T\otimes H)$ is strictly weak
mixing; \item[(ii)] $(\ga,\f,T)$ and $(\gb,\f_1,H)$ are strictly
weak mixing;
\end{itemize}
the implication (i)$\Rightarrow$(ii) holds.

If in addition one has $(\ga\otimes\gb)^*=\ga^*\otimes\gb^*$, then
(ii)$\Rightarrow$(i) also holds.
\end{thm}

\begin{pf} The implication (i)$\Rightarrow$(ii) immediately follows from the definition.

Now assume that $(\ga\otimes\gb)^*=\ga^*\otimes\gb^*$ holds. Let
us consider the implication (ii)$\Rightarrow$(i).  It is clear
that the state $\f\otimes\f_1$ is invariant with respect to
$T\otimes H$.

Let $\p\in \ga^*$ and $\v\in \gb^*$ be arbitrary functionals and
$x\in\ker\f$, $y\in\ker\f_1$. Then according to (ii) we have
\begin{equation}\label{1mx}
\lim_{n\to\infty}\frac{1}{n}\sum_{k=0}^{n-1}|\p(T^k(x))|=0.
\end{equation}

The Schwartz inequality implies that
\begin{eqnarray}\label{2mx}
\frac{1}{n}\sum_{k=0}^{n-1}|\p(T^k(x))\v(H^k(y))|&\leq&
\frac{1}{n}\sqrt{\sum_{k=0}^{n-1}|\p(T^k(x))|^2}\sqrt{\sum_{k=0}^{n-1}|\v(H^k(y))|^2}\nonumber\\
&=&
\sqrt{\frac{1}{n}\sum_{k=0}^{n-1}|\p(T^k(x))|^2}\sqrt{\frac{1}{n}\sum_{k=0}^{n-1}|\v(H^k(y))|^2}\nonumber\\
&\leq& \|\v\|\|y\|\sqrt{\frac{1}{n}\sum_{k=0}^{n-1}|\p(T^k(x))|^2}
\end{eqnarray}

Moreover, the relations
\begin{eqnarray*}
\frac{1}{n}\sum_{k=0}^{n-1}|\p(T^k(x))|^2&\leq& \sup_{0\leq k\leq
n-1}|\p(T^k(x))|\frac{1}{n}\sum_{k=0}^{n-1}|\p(T^k(x))|\\
&\leq&\|\p\|\|x\|\frac{1}{n}\sum_{k=0}^{n-1}|\p(T^k(x))|
\end{eqnarray*}
with \eqref{1mx} yield that
$$
\lim_{n\to\infty}\frac{1}{n}\sum_{k=0}^{n-1}|\p(T^k(x))\v(H^k(y))|=0.
$$
Thus,
\begin{equation}\label{eq0}
\lim_{n\to\infty}\frac{1}{n}\sum_{k=0}^{n-1}|\p\otimes\v(T^k\otimes
H^k(x\otimes y))|=0,
\end{equation}
for $x\in\ker\f$, $y\in\ker\f_1$,  $\p\in \ga^*$, $\v\in \gb^*$.

Let $\ga^*\odot\gb^*$ be the algebraic tensor product of $\ga^*$
and $\gb^*$.  Thanks to our assumption one can see that the
$\|\cdot\|_1$-closure of $\ga^*\odot\gb^*$ is $(\ga\otimes\gb)^*$.
So, using the norm-denseness of the elements
$\sum_{i=1}^m\p_i\otimes\v_i$ in $(\ga\otimes\gb)^*$ from
\eqref{eq0} one gets
\begin{equation}\label{eq1}
\lim_{n\to\infty}\frac{1}{n}\sum_{k=0}^{n-1}|\w(T^k\otimes
H^k(x\otimes y))|=0,
\end{equation}
for $\w\in (\ga\otimes\gb)^*$.

Let $x\in\ga$ and $y\in\gb$. Denoting $x^0=x-\f(x)\id$,
$y^0=y-\f_1(y)\id$ we have $x^0\in\ker\f$, $y^0\in\ker\f_1$, so
for them \eqref{eq1} holds.

Denote $\w_1(x)=\w(x\otimes\id),x\in\ga$ and $
\w_2(y)=\w(\id\otimes y)$, $y\in\gb$. Then according to condition
(ii) we find
\begin{eqnarray}\label{eq12}
\lim_{n\to\infty}\frac{1}{n}\sum_{k=0}^{n-1}|\w_1(T^k(x))-\w(\id\otimes\id)\f(x)|=0, \\
\label{eq13}
\lim_{n\to\infty}\frac{1}{n}\sum_{k=0}^{n-1}|\w_2(H^k(y))-\w(\id\otimes\id)\f_1(y)|=0.
\end{eqnarray}

Now from
\begin{eqnarray*}
\frac{1}{n}\sum_{k=0}^{n-1}|\w(T^k\otimes H^k(x\otimes
y))&-&\w(\id\otimes\id)\f(x)\f_1(y)| \nonumber \\ &\leq&
|\f_1(y)|\left(\frac{1}{n}\sum_{k=0}^{n-1}|\w_1(T^k(x))-\w(\id\otimes\id)\f(x)|\right)\nonumber
\\
&&+|\f(x)|\left(\frac{1}{n}\sum_{k=0}^{n-1}|\w_2(H^k(y))-\w(\id\otimes\id)\f_1(y)|\right)\nonumber\\
&&+\frac{1}{n}\sum_{k=0}^{n-1}|\w(T^k\otimes H^k(x^0\otimes y^0))|
\end{eqnarray*}
and \eqref{eq1}-\eqref{eq13} we get
\begin{eqnarray}\label{eq21}
\lim_{n\to\infty}\frac{1}{n}\sum_{k=0}^{n-1}|\w(T^k\otimes
H^k(x\otimes y))&-&\w(\id\otimes\id)\f(x)\f_1(y)|=0.
\end{eqnarray}

The norm-denseness of the elements $\sum_{i=1}^m x_i\otimes y_i$
in $\ga\otimes\gb$ with \eqref{eq21} yields
\begin{eqnarray*}
\lim_{n\to\infty}\frac{1}{n}\sum_{k=0}^{n-1}|\w(T^k\otimes
H^k({\mathbf{z}}))&-&\w(\id\otimes\id)\f\otimes\f_1({\mathbf{z}})|=0.
\end{eqnarray*}
for arbitrary ${\mathbf{z}}\in\ga\otimes\gb$. So,
$(\ga\otimes\gb,\f\otimes\f_1,T\otimes H)$ is strictly weak
mixing.
\end{pf}

\begin{rem} Note that analogous results for weak mixing
dynamical system defined on von Neumann algebras were proved in
\cite{L},\cite{W}.
\end{rem}

From the proved theorem we get the following

\begin{cor}
\label{mix} Let $(\ga,\f,T)$ be a state preserving $C^*$-dynamical
systems.  For the following assertions
\begin{itemize}
\item[(i)] The state preserving  $C^*$-dynamical system
$(\ga\otimes\ga,\f\otimes\f,T\otimes T)$ is uniquely ergodic;
\item[(ii)] $(\ga,\f,T)$ is strictly weak mixing;
\end{itemize}
the implication (i)$\Rightarrow$(ii) holds.

If, in addition, $\ga^*\otimes\ga^*=(\ga\otimes\ga)^*$ is
satisfied then both (i),(ii) assertions are equivalent to
\begin{itemize}
 \item[(iii)] The
state preserving  $C^*$-dynamical system
$(\ga\otimes\ga,\f\otimes\f,T\otimes T)$ is strictly weak mixing;
\end{itemize}
\end{cor}

\begin{pf} (i)$\Rightarrow$(ii). Let
$(\ga\otimes\ga,\f\otimes\f,T\otimes T)$ be uniquely ergodic. Let
$x\in\ker\f, x=x^*$. The  unique ergodicity of the dynamical
system (see Theorem \ref{ue}) implies
$$
\lim_{n\to\infty}\left\|\frac{1}{n}\sum_{k=0}^{n-1}T^k\otimes
T^k(x\otimes x)\right\|=0.
$$

Hence,
$$
\lim_{n\to\infty}\left|\frac{1}{n}\sum_{k=0}^{n-1}\p\otimes\p(T^k\otimes
T^k(x\otimes x))\right|=0,  \ \ \ \textrm{for all} \ \
\p\in\ga^*_{1,h}.
$$
Self-adjointness of $x$ yields
\begin{equation}\label{ue2}
\lim_{n\to\infty}\frac{1}{n}\sum_{k=0}^{n-1}|\p(T^k(x))|^2=0 \ \ \
\textrm{for all} \ \ \p\in\ga^*_{1,h}.
\end{equation}

By the Schwartz inequality one finds
\begin{eqnarray*}
\frac{1}{n}\sum_{k=0}^{n-1}|\p(T^k(x))|&\leq&
\frac{1}{n}\sqrt{\sum_{k=0}^{n-1}1}\sqrt{\sum_{k=0}^{n-1}|\p(T^k(x))|^2}\nonumber\\
&=& \sqrt{\frac{1}{n}\sum_{k=0}^{n-1}|\p(T^k(x))|^2},
\end{eqnarray*}
which with \eqref{ue2} implies
\begin{eqnarray}\label{eq3}
\lim_{n\to\infty}\frac{1}{n}\sum_{k=0}^{n-1}|\p(T^k(x))|=0 \ \ \ \
\textrm{for all} \ \ \forall\p\in\ga^*_{1,h}.
\end{eqnarray}

Now let $x\in\ker\f$ and $\p\in\ga^*_1$ be arbitrary. Then they
can be represented as $x=x_1+ix_2$, $\p=\p_1+i\p_2$, where
$x_1,x_2\in\ker\f$, $x_j^*=x_j$, $\p_j\in\ga^*_{1,h}$, $j=1,2$.
From
\begin{eqnarray*}
\frac{1}{n}\sum_{k=0}^{n-1}|\p(T^k(x))|\leq
\frac{1}{n}\sum_{i,j=1}^2\sum_{k=0}^{n-1}|\p_i(T^k(x_j))|
\end{eqnarray*}
and \eqref{eq3} it follows that
\begin{eqnarray}\label{eq4}
\lim_{n\to\infty}\frac{1}{n}\sum_{k=0}^{n-1}|\p(T^k(x))|=0,
\end{eqnarray}
for $x\in\ker\f$, $\p\in\ga^*_1$.

Finally let $x\in\ga$. Then the last relation \eqref{eq4} for the
element $x^0=x-\f(x)\id$ implies the assertion.

Let  $\ga^*\otimes\ga^*=(\ga\otimes\ga)^*$ be satisfied, then the
implication (ii)$\Rightarrow$(iii) is a direct consequence of
Theorem \ref{mix-a}. The implication (iii)$\Rightarrow$(i)
immediately follows from Proposition \ref{uesm}.
\end{pf}

An idea of the proof of Theorem \ref{mix-a} allows us to get some
adaptation of a result of \cite{ALW} for strictly weak mixing
dynamical systems. Namely we have the following

\begin{thm}
\label{mix-c} Let $(\ga,\f,T)$ be a state preserving
$C^*$-dynamical systems. For  the following assertions
\begin{itemize}
\item[(i)] For every state preserving uniquely ergodic
$C^*$-dynamical system $(\gb,\f_1,H)$ the state preserving
$C^*$-dynamical system $(\ga\otimes\gb,\f\otimes\f_1,T\otimes H)$
is uniquely ergodic; \item[(ii)] $(\ga,\f,T)$ is strictly weak
mixing ; \item[(iii)] For every state preserving uniquely ergodic
$C^*$-dynamical system $(\gb,\f_1,H)$ such that
$\ga^*\otimes\gb^*=(\ga\otimes\gb)^*$ the state preserving
$C^*$-dynamical system $(\ga\otimes\gb,\f\otimes\f_1,T\otimes H)$
is uniquely ergodic;
\end{itemize}
the following implications hold
(i)$\Rightarrow$(ii)$\Rightarrow$(iii).

\end{thm}

\begin{pf} (i)$\Rightarrow$(ii). According
to the condition $(\ga\otimes\gb,\f\otimes\f_1,T\otimes H)$ is
uniquely ergodic, this means that the state $\f\otimes\f_1$ is a
unique for it. Take arbitrary functional $\p\in\ga^*$ and $\v\in
\cs(\gb)$, then the unique ergodicity due to Theorem \ref{ue}
implies that
\begin{eqnarray*}
0&=&\lim_{n\to\infty}\frac{1}{n}\sum_{k=0}^{n-1}\big(\p\otimes\v(T^k\otimes
H^k(x\otimes
\id))-\p(\id)\f(x)\big)\\&=&\lim_{n\to\infty}\frac{1}{n}\sum_{k=0}^{n-1}(\p(T^k(x))-\p(\id)\f(x)).
\end{eqnarray*}
This shows unique ergodicity of $(\ga,\f,T)$. Then condition (i)
implies that $T\otimes T$ is also uniquely ergodic, therefore
Corollary \ref{mix} yields that $T$ is strictly weak mixing.

(ii)$\Rightarrow$(iii). Let $(\gb,\f_1,H)$ be a completely
positive, uniquely ergodic dynamical system such that
$\ga^*\otimes\gb^*=(\ga\otimes\gb)^*$. Then it is clear that the
state $\f\otimes\f_1$ is invariant with respect to $T\otimes H$.
The argument used in the proof of Theorem \ref{mix-a} implies that
\begin{equation}\label{c}
\lim_{n\to\infty}\frac{1}{n}\sum_{k=0}^{n-1}|\w(T^k\otimes
H^k(x\otimes y))|=0.
\end{equation}
for every $x\in\ker\f$, $y\in\gb$ and  $\w\in (\ga\otimes\gb)^*$.

Let $x\in\ga$. Then \eqref{c} holds for $x^0=x-\f(x)\id$. The
unique ergodicity of $(\gb,\f_1,H)$ implies
\begin{eqnarray}\label{eq12c}
\lim_{n\to\infty}\bigg|\frac{1}{n}\sum_{k=0}^{n-1}(\w_2(H^k(y))-\w(\id\otimes\id)\f_1(y))\bigg|=0,
\end{eqnarray}
here as before $\w_2(x)=\w(\id\otimes x)$.

Now from
\begin{eqnarray*}
\bigg|\frac{1}{n}\sum_{k=0}^{n-1}\w(T^k\otimes H^k(x\otimes
y))&-&\w(\id\otimes\id)\f(x)\f_1(y)\bigg|\nonumber \\ &\leq&
\bigg|\frac{1}{n}\sum_{k=0}^{n-1}(\w(T^k\otimes H^k)(x^0\otimes
y)\bigg|\nonumber
\\
&&+|\f(x)|\bigg|\frac{1}{n}\sum_{k=0}^{n-1}(\w_2(H^k(y))-\w(\id\otimes\id)\f_1(y))\bigg|
\end{eqnarray*}
and \eqref{c}, \eqref{eq12c} it follows that
\begin{eqnarray}\label{eq2c}
\lim_{n\to\infty}\bigg|\frac{1}{n}\sum_{k=0}^{n-1}\w(T^k\otimes
H^k(x\otimes y))-\w(\id\otimes\id)\f(x)\f_1(y)\bigg|=0
\end{eqnarray}

The density argument used  in the proof of Theorem \ref{mix-a} and
 Theorem \ref{ue} yield the required assertion.
\end{pf}

\begin{rem}\label{Ta} If in condition (i) of Theorem \ref{mix-c} we take
not all state preserving uniquely ergodic $C^*$-dynamical systems,
then the assertion of the theorem fails. Indeed, let us consider
the following example. Let $S^1=\{z\in\bc: |z|=1\}$ and $\l$ be
the Lebesgue measure on $S^1$ such that $\l(S^1)=1$. The measure
induces a positive  linear functional
$\f_\l(f)=\int\limits_{S^1}f(z)d\l(z)$ such that $\f_\l(\id)=1$.
Consider a $C^*$-algebra $\ga=C(S^1)$, where $C(S^1)$ is the space
of all continuous functions on $S^1$.  Fix an element $a=\exp(2\pi
i\a)$, where $\a\in[0,1)$ is an irrational number. Define a
mapping $T_\a:C(S^1)\mapsto C(S^1)$ by $(T_\a(f)(z))=f(az)$ for
all $f\in C(S^1)$. It is clear that $(C(S^1),\f_\l,T_\a)$ is a
state preserving $C^*$-dynamical system. Since $\a$ is irrational,
then Theorem 2 of  Chapter 3 of \cite{KSF} implies that the
defined dynamical system is uniquely ergodic. According to that
theorem the tensor product $T_\a\otimes T_\b$, acting on
$C(S^1)\otimes C(S^1)$, is also uniquely ergodic for every $\b$
which is rationally independent with $\a$.

But $T_\a$ is not strictly weak mixing. Indeed, take a linear
functional $h\in C(S^1)^*$ defined by
$h(f)=\int\limits_{S^1}zf(z)d\l(z)$, $f\in C(S^1)$. Then we have
$h(T_\a(f))=a^{-2}h(f)$ for all $f\in C(S^1)$. Thus, Proposition
\ref{mix2} implies that $T_\a$ is not strictly weak mixing.
According to Corollary \ref{mix} the tensor product $T_\a\otimes
T_\a$, acting on $C(S^1)\otimes C(S^1)$, is not uniquely ergodic.
Moreover, $T_\a\otimes T_\a$ is not ergodic. Indeed, using the
well known equality $C(S^1\times S^1)=C(S^1)\otimes C(S^1)$ we see
that $T_\a\otimes T_\a$ acts as follows
$$
(T_\a\otimes T_\a)f(x,y)=f(ax,ay), \ \ x,y\in S^1,
$$
where $f\in C(S^1\times S^1)$. For the element $g$ of $C(S^1\times
S^1)$ defined by $ g(x,y)=x/y$, we have $(T_\a\otimes
T_\a)(g)(x,y)=g(x,y)$ which means that $T_\a\otimes T_\a$ is not
ergodic.
\end{rem}

\section{Examples}

In this section we are going to provide certain examples of
strictly weak mixing ucp maps.

{\bf 1.} Let $\ga=M_2(\bc)$ and $\t$ be the normalized trace on
$\ga$. By $e_{ij}, i,j=1,2$ we denote the matrix units (in the
standard basis of $\bc$) of $\ga$. Consider
$\ce:\ga\otimes\ga\to\ga$ - the canonical conditional expectation,
i.e. $\ce(x\otimes y)=\t(y)x$ (see \cite{T}). Take
$V\in\ga\otimes\ga$ such that $\ce(VV^*)=\id$. Define
$T_V:\ga\to\ga$ by $ T_V(x)=\ce(V(\id\otimes x)V^*), \ \ x\in\ga.$
Then it is clear that $T_V$ is a ucp map with $\t(T_Vx)=\t(x)$ for
all $x\in\ga$. If its peripheral spectrum is $\{1\}$, then
$T^n_V\to\t\id$ as $n\to\infty$. In this case $(\ga,\t,T_V)$ would
be strictly weak mixing. In particular, if we choose $V$ as
follows
$$
V_\b=\sqrt{\frac{2}{1+\cosh(2\b)}}\exp\{\b(e_{12}\otimes
e_{21}+e_{21}\otimes e_{12})\}, \ \ \b\in\br
$$
then all the required conditions are satisfied.

{\bf 2.} Let $(C(K),\n,T)$ be a commutative strictly weak mixing
dynamical system. Now with the aid of  above Example 1 and Theorem
\ref{mix-a} one finds that $(C(K)\otimes
M_2(\bc),\n\otimes\t,T\otimes T_{V_\b})$ is a non-commutative
strictly weak mixing dynamical system.

{\bf 3.} First we formulate a result relating to adaptation of the
Blum-Hanson theorem (see \cite{BLRT,BH,JL,Z}) for strictly weak
mixing dynamical systems, which will be used below.

\begin{thm}\label{bhsw}
\label{mix-w} A state preserving $C^*$-dynamical system
$(\ga,\f,T)$ is strictly weak mixing if and only if
\begin{equation}\label{stmix31}
\lim_{n\to\infty}\frac{1}{n}\sum_{m=0}^{n-1}T^{k_m}(x)=\f(x)\id
\end{equation}
for every $x\in\ga$ and increasing sequence of positive numbers
$\{k_n\}$ such that $\sup_n k_n/n<\infty$. Here the convergence is
meant with respect to the uniform norm.
\end{thm}

Now let $\bbf_{\infty}$ be the free group on infinitely many
generators $\{g_{i}\}_{i\in\bz}$. Let $\l$ be the regular
representation of $\bbf_\infty$ on $\ell^2(\bbf_\infty)$. If
$\delta_t$, $t\in \bbf_\infty$ denotes the unit vectors
$$
\d_t(s)= \left\{
\begin{array}{ll}
1, \ \ s=t\\
0, \ \ s\neq t
\end{array}
\right.
$$
in $\ell^2(\bbf_\infty)$, then one has $\l(s)\d_t=\d_{st}$, for
every $s,t\in\bbf_\infty$. The $C^*$-algebra $C_\l^*(\bbf_\infty)$
associated with the regular representation of $\bbf_\infty$, is
the norm-closure in $B(\ell^2(\bbf_\infty))$ of $span\{\l(s):s\in
\bbf_\infty\}$. Note that any element $s\in\bbf_\infty$ has a
unique expression as a finite product of $g_i$ ($i\in\bz$). This
expression is called the {\it word} for $s$. The number of factors
in the word is called the {\it length} of the word. Let
$\b:\bbf_\infty\to \bbf_\infty$ be the shift-automorphism, i.e.
$\b(g_i)=g_{i+1}$ for all $i\in\bz$. The induced by $\b$
free-shift automorphism of $C_\l^*(\bbf_\infty)$ is denoted by
$\a_\b$. In \cite{AD} it has been proved that $\a_\b$ is uniformly
ergodic. Now we are going to show that it is strictly weak mixing.

By a standard density argument, it is enough to show that the
sequence $\{\a^{n}_\b(\l(s))\}_{n\geq1}$ is weakly mixing to zero
whenever $\b(s)\neq s$, that is
\begin{equation}\label{alphamix}
\frac{1}{n}\sum_{k=1}^{n}\big|f(\a^{k}_\b(\l(s)))\big|\to 0
\end{equation}
for each $f\in C^{*}_{\l}(\bbf_\infty)^{*}$.

Let $s$ be a nontrivial element of word length $p$, then by
Haagerup's inequality (cf. \cite{H}), for each sequence
$\{k_{j}\}$ of natural numbers, one has
\begin{eqnarray*}
\bigg\|\frac{1}{n}\sum_{j=1}^{n}\a^{k_{j}}_\b(\l(s))\bigg\|& \leq&
(p+1)\bigg\|\frac{1}{n}\sum_{j=1}^{n}\d_{\b^{k_{j}}(s)}\bigg\|_{\ell^{2}(\bbf_\infty)}\\
&=&\frac{p+1}{\sqrt{n}}.
\end{eqnarray*}

Now according to Theorem \ref{bhsw} we get \eqref{alphamix}.

\begin{rem} Note that in \cite{FM} some examples of strictly
weak mixing dynamical systems, related to free shift of the
reduced $C^*$-algebras of RG-groups and amalgamated free product
$C^*$-algebras have been provided.
\end{rem}

\section{Uniform weighted ergodic theorem}

From Theorem \ref{bhsw} we know that subsequential ergodic theorem
holds for strictly weak mixing dynamical system. But it would be
interesting to obtain some weighted uniform ergodic theorems. Note
that similar problem has been investigated in \cite{BLRT} for
Hilbert spaces. Namely, they found the necessary and sufficient
conditions for the convergence of
\begin{equation}\label{hw}
\frac{1}{n}\sum_{k=0}^{n-1}a_kT^kx
\end{equation}
for every contraction $T$ on a Hilbert space $H$ and every $x\in
H$. In our case, a situation is different, since we are dealing
with $C^*$-algebras, which are not Hilbert spaces. In this section
we are going to give a sufficient condition for the uniform
convergence of weighted averages \eqref{hw} for strictly weak
mixing $C^*$-dynamical systems.

By analogy of a Besicovitch sequences (see \cite{JLO}) we
introduce a notion of $S$-Besicovitch sequences as follows: we say
that a bounded sequence $\{b_n\}\subset\bc$ is {\it a
$S$-Besicovitch} if for any $\e>0$ there exists a uniquely ergodic
dynamical system $(C(K),\n,T_1)$, a function $f_0\in C(K)$ and
$\w_0\in K$ such that
\begin{equation}\label{bs}
\limsup_{n\to
\infty}\frac{1}{n}\sum_{k=0}^{n-1}\bigg|b_k-(T^k_1f_0)(\w_0)\bigg|<\e.
\end{equation}

Now we are ready to formulate the result.

\begin{thm}
\label{mix-bes} Let $(\ga,\f,T)$ be a strictly weak mixing
$C^*$-dynamical system. Then for every $x\in\ga$ and
$S$-Besicovitch sequence $\{b_n\}$ the averages
\begin{equation}\label{bs1}
\frac{1}{n}\sum_{k=0}^{n-1}b_kT^{k}(x)
\end{equation}
converge uniformly in $\ga$.
\end{thm}

\begin{pf} Let $\e>0$ be an arbitrary number. Assume that  $(C(K),\n,T_1)$, $f_0$, $\w_0$ is a generating
system for the sequence $\{b_n\}$. Due to commutativity of $C(K)$
one has  $(C(K)\otimes\ga)^*=C(K)^*\otimes\ga^*$, therefore,
Theorem \ref{mix-c} implies that a dynamical system
$(C(K)\otimes\ga,\n\otimes\f,T_1\otimes T)$ is uniquely ergodic,
i.e. for every ${\mathbf{x}}\in C(K)\otimes\ga$ the following
holds
\begin{equation*}
\lim_{n\to\infty}\frac{1}{n}\sum_{k=0}^{n-1}(T_1^k\otimes
T^{k})({\mathbf{x}})=(\n\otimes\f)({\mathbf{x}})\id.
\end{equation*}

In particular, for  $f_0\otimes x\in C(K)\otimes\ga$ we have
\begin{equation*}
\lim_{n\to\infty}\frac{1}{n}\sum_{k=0}^{n-1}(T_1^kf_0)(\w_0)T^{k}(x)=\n(f_0)\f(x)\id.
\end{equation*}

This means that there is $N_0\in\bn$ such that
\begin{equation}\label{bs2}
\bigg\|\frac{1}{n}\sum_{k=0}^{n-1}(T_1^kf_0)(\w_0)T^{k}(x)-
\frac{1}{m}\sum_{l=0}^{m-1}(T_1^lf_0)(\w_0)T^{l}(x)\bigg\|<\e
\end{equation}
for all $n,m\geq N_0$.

Now from \eqref{bs} and \eqref{bs2} we find
\begin{eqnarray*}
\bigg\|\frac{1}{n}\sum_{k=0}^{n-1}b_kT^{k}(x)-\frac{1}{m}\sum_{l=0}^{m-1}b_lT^{l}(x)\bigg\|&\leq&
\bigg\|\frac{1}{n}\sum_{k=0}^{n-1}b_kT^{k}(x)-\frac{1}{n}\sum_{k=0}^{n-1}(T_1^kf_0)(\w_0)T^{k}(x)\bigg\|\\
&&+\bigg\|\frac{1}{m}\sum_{k=0}^{m-1}b_lT^{l}(x)-\frac{1}{m}\sum_{l=0}^{m-1}(T_1^lf_0)(\w_0)T^{l}(x)\bigg\|\\
&&+\bigg\|\frac{1}{n}\sum_{k=0}^{n-1}(T_1^kf_0)(\w_0)T^{k}(x)-
\frac{1}{n}\sum_{l=0}^{m-1}(T_1^lf_0)(\w_0)T^{l}(x)\bigg\|\\
&\leq&\frac{1}{n}\sum_{k=0}^{n-1}|b_k- (T_1^kf_0)(\w_0)|\|x\|\\
&&+\frac{1}{m}\sum_{l=0}^{m-1}|b_l-
(T_1^lf_0)(\w_0)|\|x\|+\e\\
&\leq&\e(2\|x\|+1)
\end{eqnarray*}
for all $n,m\geq N_0$.
 This completes the proof.
\end{pf}

{\bf Example.} Consider the uniquely ergodic $C^*$-dynamical
system $(C(S^1),T_\a)$ defined in Remark \ref{Ta}. For fixed
$m\in\bn$ take $f_{0,m}(z)=z^m$ and $\w_0=1$. Then one can see
that a sequence $\{b^{(m)}_{n}\}_{n\in\bn}$ given by
$b^{(m)}_{n}=a^{nm}$, here as before $a=\exp\{2\pi i\a\}$, is
$S$-Besicovitch. From the just proved theorem due to
$\f_\l(f_{0,m})=0$ one concludes that for every $x\in\ga$
\begin{equation*}
\lim_{n\to\infty}\frac{1}{n}\sum_{k=0}^{n-1}a^{km}T^{k}(x)=0,
\end{equation*}
for every strictly weak mixing $C^*$-dynamical system
$(\ga,\f,T)$.

\begin{rem} We note that Besicovitch weighted ergodic type theorems were
studied in (\cite{BLRT},\cite{JLO},\cite{S}).
\end{rem}

\end{document}